\documentclass[reqno]{amsart}
\usepackage[utf8]{inputenc}
\usepackage{amsfonts, amsmath, amssymb, amsthm, graphicx, enumitem}
\usepackage[left=2cm, right=2cm, top=3cm]{geometry}
\newtheorem{theorem}{Theorem}

\newtheorem{definition}[theorem]{Definition}

\newtheorem{lemma}[theorem]{Lemma}
\newtheorem{proposition}[theorem]{Proposition}
\newtheorem{remark}[theorem]{Remark}
\usepackage{dsfont}

\renewcommand{\d}{\displaystyle}
\usepackage{graphicx}
\usepackage{apptools}
\AtAppendix{\counterwithin{lem}{section}}
\newtheorem{lem}{Lemma}
\newcommand{\pts}[1]{\left(#1\right)}  
\newcommand{\cts}[1]{\left[#1\right]}                                  	  
\newcommand{\lvs}[1]{\left\{#1\right\}}                                	  
\newcommand{\abs}[1]{\left|#1\right|}                                  	  

\newcommand{\C}{\mathbb{C}} 
\newcommand{\R}{\mathbb{R}}
\newcommand{\A}{\mathcal{A}}
\newcommand{\al}{\alpha}

\newcommand{\U}{\mathcal{U}}
\newcommand{\OO}{\mathcal{O}}
\newcommand{\Fi}{\Phi}
\newcommand{\D}{C_c^\infty(0,1)}
\newcommand{\HH}{\mathcal{H}}

\begin{document}
\title[Null control of a degenerate equation with drift]{Boundary controllability for a 1D degenerate parabolic equation with drift and a singular potential
}

\author[L. Galo-Mendoza \and M. L\'opez-Garc\'ia]{Leandro Galo-Mendoza \and Marcos L\'opez-Garc\'ia}
\address{Instituto de Matem\'{a}ticas-Unidad Cuernavaca\\
Universidad Nacional Aut\'onoma de M\'exico\\
Av. Universidad S/N\\
Cuernavaca, Morelos, C.P. 62210\\
M\'{e}xico.}
\email{jesus.galo@im.unam.mx, marcos.lopez@im.unam.mx}

\maketitle

\begin{abstract}We prove the null controllability of a one dimensional degenerate parabolic equation with drift and a singular potential. We study the case the potential arises at the left end point and the weighted Dirichlet boundary control is located at this point. We get a spectral decomposition of a suitable operator, defined in a weighted Sobolev space, involving Bessel functions and their zeros, then we use the moment method by Fattorini and Russell to obtain an upper estimate of the cost of controllability. We also obtain a lower estimate of the cost of controllability by using a representation theorem for analytic functions of exponential type.
\end{abstract}

\section{Introduction and main results}
Let $T > 0$ and set $Q := (0, 1) \times (0, T )$. For $\al,\beta\in \R$ with $0\leq\al<2$, $\al+\beta <1$, consider the system
\begin{equation}\label{problem}
\left\{\begin{aligned}
u_t-(x^\al u_x)_x-\beta x^{\al -1}u_x-\frac{\mu}{x^{2-\al}} u&=0 & & \text { in }Q, \\
\pts{x^{-\gamma}u}(0, t) =f(t), u(1, t)&=0 & & \text { on }(0, T), \\
u(x, 0) &=u_{0}(x) & & \text { in }(0, 1),
\end{aligned}\right.
\end{equation}
provided that $\mu\in\R$ satisfies
\begin{equation}\label{mucon}
-\infty<\mu<\mu(\al+\beta),
\end{equation}
where
\begin{equation}\label{gamadef}
\mu(\delta):=\frac{(1-\delta)^2}{4}, \quad \delta\in\R,\quad \text{and}\quad\gamma=\gamma(\alpha, \beta,\mu):=\sqrt{\mu(\alpha+\beta)}-\sqrt{\mu(\alpha+\beta)-\mu}.
\end{equation}

The aim of this work is twofold. First, we provide a notion of weak solution for system (\ref{problem}) and establish the well-posedness of this problem in suitable interpolation spaces. Notice we consider a weighted Dirichlet boundary condition at the left end point in order to compensate the singularity of the potencial at this point. Second, we are interested in the null controllability of the system (\ref{problem}). Here we apply the moment method by Fattorini and Russell \cite{Fatorini} to prove the null controllability and show an upper bound estimate of the cost of controllability. At the end, we use a representation theorem for analytic functions of exponential type to obtain a lower bound estimate of the cost of controllability.\\

This paper generalizes several results about null controllability for 1D degenerate/singular parabolic equations when the control acts at the end point where the degeneracy/singularity arises. In this setting, the first result was obtained in \cite{gueye} in the case of a weak degeneracy $0 \leq \alpha< 1$, without singularity and drift, i.e $\beta=\mu=0$. In that work the null controllability is achieved by means of $L^2$ controls.\\

The last result was complemented in \cite{canna} where sharp estimates for the cost of the control were obtained. In \cite{bic} the author analyzes the controllability of the heat equation with a singular inverse-square potential, i.e with $\alpha=\beta=0$ and $\mu<\mu(0)=1/4$; In \cite{bic2} the authors also study the weakly degenerate case without drift ($\beta=0$) and $\mu$ satisfying (\ref{mucon}). All these results get the null controllability by means of $H^1$ controls. In \cite{bic,bic2} the authors justify the choice of the weighted Dirichlet boundary condition in system (\ref{problem}). Unfortunately, the proofs of \cite[Theorem 2.1]{bic} and \cite[Theorem 2.2]{bic2} are not complete, see Remark \ref{error} for details.\\

With respect to the strongly degenerate case ($1<\alpha<2$) with no singularity ($\mu=0$), in \cite{gueye2}  the authors analyze the null controllability of a degenerate parabolic equation with a degenerate one-order transport term. This work motivate us to add a first order term along with a suitable drift in a degenerate/singular parabolic equation. In this case, the first order term allow us to consider both cases, weak and strong degeneracy, without changing the (weighted) Dirichlet boundary condition. In \cite{du,flores,flores2} the authors consider the null controllability of 1D degenerate parabolic equations with first order terms, but they use interior controls.\\

Now assume system (\ref{problem}) admits a unique solution for initial conditions in a Hilbert space $H$, that it is described in the next section. It is said that system (\ref{problem}) is null controllable in $H$ at time $T>0$ with controls in $L^2(0,T)$ if for any $u_0 \in H$ there exists $f \in L^2(0,T)$ such that the corresponding solution satisfies $u(\cdot,T ) \equiv 0$.\\

Once we know system (\ref{problem}) is null controllable we are interested in the behavior of the cost of the controllability. Thus, consider the set of admissible controls
\[U(T,\alpha,\beta,\mu,u_0)=\{f \in L^2(0, T ): u \text{ solution of system (\ref{problem}) that satisfies } u(\cdot,T ) \equiv 0\}.\]
Then the cost of the controllability is defined as 
\[\mathcal{K}(T,\alpha,\beta,\mu):=\sup_{\|u_0\|_H\leq1}\inf\left\{\|f\|_{L^2(0,T)}:f\in U(T,\alpha,\beta,\mu,u_0)\right\}.\]

The main result of this work is the following.
\begin{theorem}
Let $T>0$ and $\al,\beta,\mu,\gamma\in \R$ with $0\leq\al<2$, $\al+\beta <1$, $\mu$ and $\gamma$ satisfying (\ref{mucon}) and (\ref{gamadef}) respectively. The next statements hold.
\begin{enumerate}
\item \textbf{Existance of a control} For any $f \in L^2(0,T)$ and $u_0\in L^2((0,1);x^\beta dx)$ there exists a control $f \in L^2(0, T )$ such that the solution $u$ to (\ref{problem}) satisfies $u(\cdot,T ) \equiv 0$.
\item \textbf{Upper bound of the cost} There exists a constant $c>0$ such that for every $\delta\in (0,1)$ we have
		\[\mathcal{K}(T,\alpha,\beta,\mu)\leq \frac{c M(T,\alpha,\nu,\delta)T^{1/2}\kappa_\al^{-1/2}}{\sqrt{\mu(\alpha+\beta)}+\sqrt{\mu(\alpha+\beta)-\mu)}}
\exp\pts{-\frac{T}{2}\kappa_\alpha^2 j_{\nu,1}^2},\]
where 
\begin{equation}\label{Nu}
\kappa_\al:=\frac{2-\al}{2},\quad \nu=\nu(\al,\beta,\mu):=\sqrt{\mu(\al+\beta)-\mu}/\kappa_\al,
\end{equation}
$j_{\nu,1}$ is the first positive zero of the Bessel function $J_\nu$ (defined in the Appendix), and
\[M(T,\alpha,\nu,\delta)=\pts{1+\frac{1}{(1-\delta)\kappa_\alpha^2 T}}\cts{\exp\pts{\frac{1}{\sqrt{2}\kappa_\alpha}}+\frac{1}{\delta^3}\exp\pts{\frac{3}{(1-\delta)\kappa_\alpha^2 T}}}\exp\pts{-\frac{(1-\delta)^{3/2}T^{3/2}}{8(1+T)^{1/2}}\kappa_\alpha^3 j_{\nu,1}^2}.\]
\item \textbf{Lower bound of the cost} There exists a constant $c>0$ such that
		\[\frac{c2^{\nu} \Gamma(\nu+1) \left|J_{\nu}^{\prime}\left(j_{\nu, 1}\right)\right|\exp{\left(\pts{\frac{1}{2}-\frac{\log 2}{\pi}}j_{\nu,2}\right)}}{\pts{{2T \kappa_\alpha}}^{1/2}\pts{\sqrt{\mu(\alpha+\beta)}+\sqrt{\mu(\alpha+\beta)-\mu}}\left(j_{\nu, 1}\right)^{\nu}}\exp\pts{-\pts{j_{\nu,1}^2+\frac{j_{\nu,2}^2}{2}}\kappa_\alpha^2 T}\leq \mathcal{K}(T,\alpha,\beta,\mu),\]
\end{enumerate}
where $j_{\nu,2}$ is the second positive zero of the Bessel function $J_\nu$.
\end{theorem}

\section{Functional setting and well-posedness}
In order to get our results we will consider suitable weighted spaces. First, we introduce the weighted Lebesgue space $L^2_\beta(0,1):=L^2((0,1);x^\beta dx)$, $\beta\in\R$, with the inner product
\[\langle f,g\rangle_\beta:=\pts{\int_0^1 f(x)g(x)x^\beta \mathrm{d}x}^{1/2},\]
and its corresponding norm denoted by $\|\cdot\|_{\beta}$.\\

For $\al,\beta\in \R$ we also introduce the weighted Sobolev space
\[H_{\alpha, \beta}^{1}(0,1)=\left\{u \in L_{\beta}^{2}(0,1)\cap H^1_{\textrm{loc}}(0,1): x^{\alpha / 2} u_{x} \in L_{\beta}^{2}(0,1)\right\}\]
with the inner product
\[\langle u, v\rangle_{\alpha, \beta}:=\int_{0}^{1} u v\, x^{\beta} \!\mathrm{d}x+\int_{0}^{1} x^{\alpha+\beta} u_{x} v_{x} \mathrm{d}x,\]
and its corresponding norm denoted by $\|\cdot\|_{\alpha,\beta}$.\\

Under suitable conditions on the parameters $\alpha, \beta$, we can see that $H_{\alpha, \beta}^{1}(0,1)$ is contained in the space of locally absolutely continuous functions on $(0,1)$, so we can talk about the trace at $x=0,1$ of functions in $H_{\alpha, \beta}^{1}(0,1)$.
\begin{remark} \label{abscont}
Let $\al,\beta\in \R$ with $\al\geq0, \al+\beta <1$. We have that $H_{\alpha, \beta}^{1}(0,1) \subset W^{1,1}(0,1)$. This follows from the inequalities
\[\int_{0}^{1}|u| \mathrm{d}x \leq \frac{1}{(1-\beta)^{1 / 2}}\left(\int_{0}^{1} |u|^{2}  x^{\beta}\mathrm{d}x\right)^{1 / 2},\quad\text{and  }\]
\[\int_{0}^{1}\left|u_{x}\right| \mathrm{d}x \leq \frac{1}{(1-\alpha-\beta)^{1 / 2}}\left(\int_{0}^{1} x^{\alpha+\beta}\left|u_{x}\right|^{2} \mathrm{d}x\right)^{1 / 2}\]
for all $u\in H_{\alpha, \beta}^{1}(0,1)$. This implies the existence of the limits $u(0):=\lim_{x\rightarrow 0^+}u(x),$ and $u(1):=\lim_{x\rightarrow 1^-}u(x)$.  In particular $H_{\alpha, \beta}^{1}(0,1)\subset C([0,1])$.
\end{remark}

\begin{definition}
For $\al,\beta\in \R$ with $\al\geq 0$, $\al+\beta <1$, we consider the space
\[H_{\alpha, \beta,0}^{1}=H_{\alpha, \beta,0}^{1}(0,1):=\left\{ u \in H_{\alpha, \beta}^{1}(0,1): u(0)=u(1)=0\right\}.\]
\end{definition}

For $u\in H_{\alpha, \beta,0}^{1}$ notice that
\begin{equation}\label{est1}
\left|u(x)\right|\leq \int_{0}^{x} |u_{x}(s)| \mathrm{d}s \leq \frac{x^{(1-\alpha-\beta) / 2}}{(1-\alpha-\beta)^{1 / 2}}\left(\int_{0}^{x} s^{\alpha+\beta}\left|u_{x}(s)\right|^{2} \mathrm{d}s\right)^{1/2},
\end{equation}
therefore
\begin{equation}\label{bnd0}
\lim _{x \rightarrow 0^{+}} x^{\delta}|u(x)|=0\quad \text{for all  }\delta \geq(\alpha+\beta-1) / 2.
\end{equation}
 
 In the setting of one-dimensional degenerate/singular operators a useful tool to consider is the so-called generalized Hardy inequality,
\begin{equation}\label{hardy}
\mu(\delta)\int_0^1\frac{u^2}{x^{2-\delta}}\mathrm{d}x\leq\int_0^1x^\delta u_x^2\mathrm{d}x
\end{equation}
for all $u\in C^\infty_c(0,1)$ (the space of infinitely smooth functions compactly supported in $(0,1)$), with $\delta\in \R.$\\

 Next, we generalize the Hardy inequality in the setting of the weighted Sobolev space $H_{\alpha, \beta,0}^{1}$.
 \begin{proposition}\label{Hardy}
For $\al,\beta\in \R$ with $\al\geq 0$, $\al+\beta <1$ the inequality (\ref{hardy}) holds for any $u\in H_{\alpha, \beta,0}^{1}$, with $\delta=\al+\beta$.
 \end{proposition}
 \begin{proof}
 Let $u\in H_{\alpha, \beta,0}^{1}$. Remark \ref{abscont} implies that $u^2\in W^{1,1}(\varepsilon,1) $ provided that $\varepsilon>0$ is small enough. Set $\delta=\al+\beta$ then
 \begin{eqnarray*}
 \int_\varepsilon^1\pts{x^{\delta/2}u_x-\frac{1-\delta}{2}\frac{u}{x^{(2-\delta)/2}}}^2\mathrm{d}x
 &=&\int_\varepsilon^1x^{\delta}|u_x|^2\mathrm{d}x+\mu(\delta)\int_\varepsilon^1\frac{u^2}{x^{2-\delta}}\mathrm{d}x-\frac{1-\delta}{2}\int_\varepsilon^1\frac{(u^2)_x}{x^{1-\delta}}\mathrm{d}x\\
 &=& \int_\varepsilon^1x^{\delta}|u_x|^2\mathrm{d}x-\mu(\delta)\int_\varepsilon^1\frac{u^2}{x^{2-\delta}}\mathrm{d}x-\frac{1-\delta}{2}\pts{\lim_{x\rightarrow 1^{-}}\frac{|u(x)|^2}{x^{1-\delta}}-\frac{|u(\varepsilon)|^2}{\varepsilon^{1-\delta}}},
 \end{eqnarray*}
 and the result follows by taking the limit as $\varepsilon \rightarrow 0^+$ and using the equality (\ref{bnd0}).
 \end{proof}
 
Now we assume $0\leq \al<2$, $\al+\beta<1$. For any $u \in H_{\alpha, \beta, 0}^{1}$, and proceeding as in (\ref{est1}), we obtain the weighted Poincaré inequality
\begin{equation}\label{poincare}
\int_{0}^{1}|u|^{2} x^{\beta} \mathrm{d}x \leq \frac{1}{(2-\alpha)(1-\alpha-\beta)} \int_{0}^{1} x^{\alpha+\beta}\left|u_{x}\right|^{2} \mathrm{d}x,
\end{equation}
therefore
\[\|u\|_{\al,\beta,0}:=\left(\int_{0}^{1} x^{\alpha+\beta}\left|u_{x}\right|^{2} \mathrm{d}x\right)^{1 / 2}\]
is an equivalent norm to $\|u\|_{\alpha, \beta}$ in $H_{\al, \beta, 0}^{1}$.\\

For $\mu<\mu(\alpha+\beta)$, Proposition \ref{Hardy} implies that
\[
\|u\|_{*}=\left(\int_{0}^{1} x^{\al+\beta}\left[\left|u_{x}\right|^{2}-\frac{\mu}{x^{2}} u^{2}\right] \mathrm{d}x\right)^{1 / 2}
\]
is an equivalent norm to $\|u\|_{\al,\beta,0}$ in $H_{\alpha, \beta, 0}^{1}$. In fact, we have
 \begin{eqnarray*}
\|u\|_{\al,\beta,0} &\leq & \| u\|_{*} \leq \left(1-\frac{\mu}{\mu(\alpha+\beta)}\right)^{1/2}\| u \|_{\al,\beta,0},\quad \mu<0, \\
\left(1-\frac{\mu}{\mu(\alpha+\beta)}\right)^{1/2}\|u\|_{\al,\beta,0} &\leq & \| u\|_{*} \leq \left\|u \right\|_{\al,\beta,0},\quad 0\leq \mu <\mu(\al+\beta).
\end{eqnarray*}

Since $\D\subset H_{\alpha, \beta, 0}^{1}\subset L^2_\beta(0,1),$ and (\ref{poincare}) implies that the inclusion $(H_{\alpha, \beta,0}^1,\|\cdot\|_*)\hookrightarrow L^2_\beta(0,1)$ is continuous, the following definition makes sense.

\begin{definition}
For $\al,\beta\in \R$ with $0\leq\al<2$, $\al+\beta <1$, we consider the Gelfand triple $\left((H_{\alpha, \beta,0}^1,\|\cdot\|_*), L^2_\beta(0,1), H_{\alpha, \beta,0}^{-1}\right)$, i.e $H_{\alpha, \beta,0}^{-1}$ stands for the dual space of $(H_{\alpha, \beta,0}^1,\|\cdot\|_*)$ with respect to the pivot space $L^2_\beta(0,1)$:
\[(H_{\alpha, \beta,0}^1,\|\cdot\|_*)\hookrightarrow L^2_\beta(0,1)=\left(L^2_\beta(0,1)\right)'\hookrightarrow H_{\alpha, \beta,0}^{-1}:=\pts{H_{\alpha, \beta,0}^1,\|\cdot\|_*}'.\]
\end{definition}

The inner product $\langle \cdot ,\cdot \rangle_*$ induces an isomorphism $\A:H_{\alpha, \beta,0}^1\rightarrow H_{\alpha, \beta,0}^{-1}$ given by
\[\langle u ,v \rangle_*=\langle \A u,v\rangle_{H_{\alpha, \beta,0}^{-1},H_{\alpha, \beta,0}^{1}},\quad u,v\in H_{\alpha, \beta,0}^{1}.\]

Let $D(\A):=\A^{-1}(L^2_\beta(0,1))=\{u\in H_{\alpha, \beta,0}^{1}: \A u\in L^2_\beta(0,1)\}=\{u\in H_{\alpha, \beta,0}^{1}: \exists f\in L^2_\beta(0,1) \text{ such that }\langle u,v \rangle_*=\langle f,v\rangle_\beta \forall v\in H_{\alpha, \beta,0}^{1} \}$.\\

\begin{proposition} For $\al,\beta,\mu\in \R$ with $0\leq\al<2$, $\al+\beta <1$ and $\mu<\mu(\alpha+\beta)$, we have
\[D(\A)=\left\{u\in H^1_{\al,\beta,0}\cap H^2_{loc}(0,1): (x^\al u_x)_x+\beta x^{\al -1}u_x+\frac{\mu}{x^{2-\al}} u\in L^2_\beta(0,1)\right\}.\]
\end{proposition}
\begin{proof}
Let $H$ be the set in the right hand side, we will show that $D(\A)=H$.\\

Pick $u \in D(\A)$, then there exists $f \in L^{2}_\beta(0,1)$ such that 
\[\int_{0}^{1} \left(x^{\alpha+\beta} u_{x} v_{x}-\frac{\mu}{x^{2-\alpha-\beta}} u v\right)\mathrm{d}x=\int_0^1 f v x^{\beta} \mathrm{d}x
\quad \text{for all }v \in H_{\alpha, \beta, 0}^{1}.
\]
In particular,
\[
\int_{0}^{1} x^{\alpha+\beta} u_{x} v_{x}\mathrm{d}x=\int_{0}^{1}\left(f+\frac{\mu}{x^{2-\alpha}} u\right) v x^{\beta} \mathrm{d}x\quad \text {for all  } v \in \D,
\]
hence
\[-\left(x^{\alpha+\beta} u_{x}\right)_{x}=\left(f+\frac{\mu}{x^{2-\alpha}}u\right) x^{\beta} \text { in } \D^{\prime},
\] 
which implies
\[
\left(x^{\alpha} u_{x}\right)_{x}+\beta x^{\alpha-1} u_{x}+\frac{\mu}{x^{2-\alpha}} u=-f \text { in } \D^{\prime},
\]
therefore $u\in H$.\\

Now let  $u \in H$. For all $v \in H_{\alpha, \beta, 0}^{1} $ we claim that $x^{\alpha+\beta} u_{x} v \in W^{1,1}(0,1)$:
\[\int_0^1x^{\al+\beta}|u_x v|\mathrm{d}x\leq \pts{\int_0^1x^{\al+\beta}|u_x|^2\mathrm{d}x}^{1/2}\pts{\int_0^1x^{\al+\beta}|v|^2\mathrm{d}x}^{1/2}\leq \|u\|_{\al,\beta,0}\|v\|_\beta,\]
\begin{equation}\label{partes}
\left(x^{\alpha+\beta} u_{x} v\right)_{x}= x^{\alpha+\beta} u_{x} v_{x}+x^{\beta}\left(\left(x^{\alpha} u_{x}\right)_x+\beta x^{\alpha-1} u_{x}+\frac{\mu}{x^{2-\alpha}} u\right)v -\frac{\mu}{x^{2-\alpha-\beta}} u v \in L^{1}(0,1),
\end{equation}
therefore there exists
\[
L=\lim _{x \rightarrow 0^{+}} x^{\alpha+\beta} u_x(x) v(x).
\]
 If $L \neq 0$ then (\ref{est1}) implies that
 \[
 x^{\alpha+\beta}\left|u_{x}(x)\right| \geq \frac{|L|}{2|v(x)|} \geq C_{\alpha, \beta} x^{(\alpha+\beta-1) / 2}
 \quad\text{as }x\rightarrow 0^+\Longrightarrow x^{\alpha+\beta}\left|u_{x}(x)\right|^2 \geq C_{\alpha, \beta} x^{-1}\quad\text{as }x\rightarrow 0^+,
 \]
 which is a contradiction to $u \in H_{\alpha, \beta, 0}^{1}$. \\
 
 On the other hand, since $u \in H$ then $u\in H^2(1/2,1)$, which implies $\lim _{x \rightarrow 1^{-}} x^{\alpha+\beta} u_x(x) v(x)=0$ for all $v \in H_{\alpha, \beta, 0}^{1} $.\\
 
  Thus, from (\ref{partes}) we have
 \begin{equation}\label{quasi}
 \int_{0}^{1} \left(x^{\alpha+\beta} u_{x} v_{x}-\frac{\mu}{x^{2-\alpha-\beta}} u v\right) \mathrm{d}x=-\int_{0}^{1} x^{\beta}\left(\left(x^{\alpha} u_{x}\right)_{x}+\beta x^{\alpha-1} u_{x}+\frac{\mu}{{x^{2-\alpha}}} u\right) v \mathrm{d}x
 \end{equation}
for all $u \in H, v \in H_{\alpha, \beta, 0}^{1}$. Therefore $u \in D(\A)$.
\end{proof}

For $\al,\beta,\mu\in \R$ with $0\leq\al<2$, $\al+\beta <1$, $\mu<\mu(\al+\beta)$, we consider the unbounded operator $\A:D(\A)\subset L^2_\beta(0,1)\rightarrow L^2_\beta(0,1)$ given by
\[\A u:=-(x^\al u_x)_x-\beta x^{\al -1}u_x-\frac{\mu}{x^{2-\al}} u.\]
From Proposition 9 in \cite[p. 370]{dautray} we have that $\A$ is a closed operator with $D(\A)$ dense in $L^2_\beta(0,1)$. We also have that $\A:(D(\A),\|\cdot\|_{D(\A)})\rightarrow L^2_\beta(0,1)$ is an isomorphism, where
\[\|u\|_{D(\A)}=\|u\|_\beta+\|\A u\|_\beta,\quad u\in D(\A).\]

Next, for any $\nu>0$ consider the operator $\A_0:D(\A_0)\subset L^2(0,1)\rightarrow L^2(0,1)$ given by
\[\A_0 u=-u_{xx}-\frac{1/4-\nu^2}{x^2}u,\quad u\in D(\A_0),\]
where
\[D(\A_0)=\left\{u\in H^1_0(0,1)\cap H^2_{\textrm{loc}}(0,1):u_{xx}+\frac{1/4-\nu^2}{x^2}u\in L^2(0,1)\right\}.\]

As far as we know the following result is new. As a consequence, we get that $\mathcal{A}$ y $\mathcal{A}_0$ share spectral properties.
\begin{proposition}
Let $\al,\beta,\mu,\nu\in \R$ with $0\leq\al<2$, $\al+\beta <1$, $\mu<\mu(\al+\beta)$, and $\nu$ given in (\ref{Nu}). Then the operator $\mathcal{A}$ is unitarily equivalent to $\kappa^2_\alpha \mathcal{A}_0$.
\end{proposition}
\begin{proof}
Consider the unitary operator $\mathcal{U}:L^2(0,1)\rightarrow L^2_\beta(0,1)$ given by
\begin{equation}\label{unitary}
\U u(x):=\kappa_\al^{1/2}x^{-\al/4-\beta/2}u(x^{\kappa_\al}), \quad u\in L^2(0,1).
\end{equation}
We claim that $u\in D(\mathcal{A}_0)$ iff $\U u\in D(\mathcal{A})$, and $\kappa^2_\alpha \mathcal{A}_0=\U^{-1}\mathcal{A}\U$.\\

Let $w:=\U u$ with $u$ a smooth function, thus 
\[u(y)=\frac{1}{\kappa_\al^{1/2}}x^{\al/4+\beta/2}w(x)\quad  \text{with  }y=x^{\kappa_\al}.  \]
By using that
\[\frac{\mathrm{d}x}{\mathrm{d}y}=\frac{1}{\kappa_\al}x^{1-\kappa_\al}=\frac{1}{\kappa_\al}x^{\al/2},\quad \frac{\mathrm{d}^2x}{\mathrm{d}y^2}=\frac{\al}{2\kappa_\al^{2}}x^{\al-1},\]
easy but tedious computations yield
\begin{equation}\label{onederi}
u_y(y)=\frac{1}{\kappa_\al^{1/2}}x^{\al/4+\beta/2}\left(\left(\frac{\al}{4}+\frac{\beta}{2}\right)\frac{w(x)}{x}+w_x(x)\right)\frac{\mathrm{d}x}{\mathrm{d}y},
	 \quad\text{and}
\end{equation}
\begin{eqnarray}
u_{yy}(y)&=&\frac{x^{\al/4+\beta/2}}{\kappa_\al^{1/2}}\left[\left(\left(\frac{\al}{4}+\frac{\beta}{2}\right)\frac{w}{x}+w_x\right)\frac{\mathrm{d}^2x}{\mathrm{d}y^2}+\left(\frac{\mathrm{d}x}{\mathrm{d}y}\right)^2\left(w_{xx}+\left(\frac{\al}{2}+\beta\right)\frac{w_x}{x} +\left(\frac{\al}{4}+\frac{\beta}{2}\right)\left(\frac{\al}{4}+\frac{\beta}{2}-1\right)\frac{w}{x^2}\right)\right]\notag \\
		&=& \frac{1}{\kappa_\al^{5/2}}x^{\al/4+\beta/2}\left[x^\al w_{xx} +(\al+\beta)x^{\al-1} w_x+\left(\frac{\al}{2}\left(\frac{\al}{4}+\frac{\beta}{2}\right)+ \left(\frac{\al}{4}+\frac{\beta}{2}\right)\left(\frac{\al}{4}+\frac{\beta}{2}-1\right) \right)x^{\al-2}w \right]\notag\\
		&=& \frac{1}{\kappa_\al^{5/2}}x^{\al/4+\beta/2}\left[x^\al w_{xx} +(\al+\beta)x^{\al-1} w_x+\left(\mu(\alpha+\beta)-\kappa_\alpha^2/4 \right)x^{\al-2}w \right],\label{twoderi}
\end{eqnarray}
because $4\mu(\alpha+\beta)=\pts{(1-\alpha/2)-(\alpha/2+\beta)}^2$.\\

Since $u\in H^2_{\textrm{loc}}(0,1)$ iff $u,u_x,u_{xx}\in L^2_{\textrm{loc}}(0,1)$, the last computations imply that $u\in H^2_{\textrm{loc}}(0,1)$ iff $w=\U u\in H^2_{\textrm{loc}}(0,1)$, also $u\in H^1_{\textrm{loc}}(0,1)$ iff $w=\U u\in H^1_{\textrm{loc}}(0,1)$.\\

Assume that $u\in H_0^1(0,1)$. Then (\ref{bnd0}) implies that $x^\delta u(x)\rightarrow 0$ as $x\rightarrow 0^+$, provided that $\delta\geq -1/2$, therefore
\[\U u(0)=\lim _{x \rightarrow 0^{+}}\U u(x)=\kappa_\al^{1/2}\lim _{y \rightarrow 0^{+}}y^{\frac{2}{\al -2}\pts{\frac{\al}{4}+\frac{\beta}{2}}}u(y)=0,\quad \text{and}\quad \U u(1)=\kappa_\al^{1/2}u(1)=0.\]
Notice that $w=\U u\in H^1_{\textrm{loc}}(0,1)\cap L^2_\beta(0,1)$. By using (\ref{onederi}), and proceeding as in the proof of Proposition \ref{Hardy}, we have
\begin{eqnarray*}
0\leq\kappa_\al^2\int_0^1|u_y(y)|^2\mathrm{d}y&=&\lim_{\varepsilon\rightarrow 0^+}\int_\varepsilon^{1-\varepsilon}x^{\al+\beta}\pts{|w_x|^2+\left(\frac{\al}{4}+\frac{\beta}{2}\right)^2\frac{|w|^2}{x^2}+\left(\frac{\al}{4}+\frac{\beta}{2}\right)\frac{(w^2)_x}{x}}\mathrm{d}x\\
						&=&\lim_{\varepsilon\rightarrow 0^+}\int_\varepsilon^{1-\varepsilon}x^{\al+\beta}\pts{|w_x|^2+\pts{\left(\frac{\al}{4}+\frac{\beta}{2}\right)^2-(\al+\beta-1)\left(\frac{\al}{4}+\frac{\beta}{2}\right)}\frac{|w|^2}{x^2}}\mathrm{d}x\\
						&&+\left(\frac{\al}{4}+\frac{\beta}{2}\right)\kappa_\al \lim_{\varepsilon\rightarrow 0^+}\left.\frac{x^{-\al/2-\beta}u^2(x^{\kappa_\al})}{x^{1-\al-\beta}}\right|_{\varepsilon}^{1-\varepsilon}\\
						&=&\int_0^{1}x^{\al+\beta}\pts{|w_x|^2+\pts{\left(\frac{\al}{4}+\frac{\beta}{2}\right)^2-(\al+\beta-1)\left(\frac{\al}{4}+\frac{\beta}{2}\right)}\frac{|w|^2}{x^2}}\mathrm{d}x<\infty.
\end{eqnarray*}
From Proposition \ref{Hardy} we get
\[\int_0^{1} \frac{|w(x)|^2}{x^{2-\al-\beta}}\mathrm{d}x= \int_0^{1} \frac{|u(y)|^2}{y^{2}}\mathrm{d}y\leq 4\int_0^{1} |u_y(y)|^2\mathrm{d}y<\infty.\]
The last two inequalities imply that $x^{\al/2}w_x\in L^2_\beta(0,1)$. Therefore $w=\U u\in H^1_{\al,\beta,0}$.\\

Conversely, assume that $w=\U u\in H^1_{\al,\beta,0}$. Then $u\in L^2(0,1)\cap H^1_{\textrm{loc}}(0,1)$. Since $x^{\al/2}w_x\in L^2_\beta(0,1)$, we can apply Proposition \ref{Hardy} and use (\ref{onederi}) to obtain that $u_y\in L^2(0,1)$. Furthermore, $u(1)=\kappa_\al^{-1/2}w(1)=0$, and $u(0)=\lim_{y\rightarrow 0^+}u(y)=\kappa_\al^{-1/2}\lim_{x\rightarrow 0^+}x^{\al/4+\beta/2}w(x)=0$, by (\ref{bnd0}). Therefore $u\in H^1_0(0,1)$. \\

In conclusion, $u\in H^1_0(0,1)$ iff $\U u\in H^1_{\al,\beta,0}$. \\

Assume that $u\in H^2_{\textrm{loc}}(0,1)$ (iff $w=\U u\in H^2_{\textrm{loc}}(0,1)$). Equality (\ref{twoderi}) implies that
\[\A_0u(y)=-u_{yy}(y)-\frac{1/4-\nu^2}{y^2}u(y)=\frac{1}{\kappa_\al^{5/2}}x^{\al/4+\beta/2}\A w(x),\]
therefore $\A_0 u\in L^2(0,1)$ iff $\A w\in L^2_\beta(0,1)$; hence $u\in D(\mathcal{A}_0)$ iff $\U u\in D(\mathcal{A})$, and
\begin{equation}\label{similar}
\kappa_\al^2\U\mathcal{A}_0u=\mathcal{A}\U u.
\end{equation}
\end{proof}

The next result shows that $\mathcal{A}$ is a diagonalizable operator whose Hilbert basis of eigenfunctions can be written in terms of a Bessel function of the first kind $J_{\nu}$ and its corresponding zeros $j_{\nu,k}$, $k\geq 1$, located in the positive half line. In the appendix we give some properties of Bessel functions and their zeros.
\begin{proposition}
$-\A$ is a negative self-adjoint operator. Furthermore, the family
\[\Fi_k(x):=\frac{(2\kappa_\al)^{1/2}}{\abs{J'_\nu\pts{j_{\nu,k}}}}x^{\pts{1-\al-\beta}/2}J_\nu\pts{j_{\nu,k}x^{\kappa_\al}},\quad k\geq 1,\]
is an orthonormal basis for $L^2_\beta(0,1)$ such that
\begin{equation}\label{lambdak}
\mathcal{A}\Fi_k=\lambda_k \Fi_k, \quad \lambda_k=\kappa_\al^2 \pts{j_{\nu,k}}^2,\quad k\geq 1,
\end{equation}
where $\nu$ is defined in (\ref{Nu}).
\end{proposition}
\begin{proof} 
From (\ref{quasi}) we get that $\A$ is a symmetric operator. Letting $u=v\in D(\A)$ in (\ref{quasi}) and using Proposition \ref{Hardy} we obtain that $-\A\leq 0$.\\

We claim that $\textrm{Ran}(I+\A)=L^2_\beta(0,1)$: Let $f \in L_{\beta}^{2}(0,1)$ be given. Since the inner product $\langle\cdot,\cdot\rangle_\beta+\langle\cdot,\cdot\rangle_*$ is equivalent to $\langle\cdot,\cdot\rangle_{\alpha,\beta}$ in $H_{\alpha, \beta, 0}^{1}$ and $f\in (H_{\alpha, \beta, 0}^{1},\|\cdot\|_{\alpha,\beta})'$, the Riesz representation theorem implies that there exists a unique $u\in H_{\alpha, \beta, 0}^{1}$ such that
\[\int_{0}^{1} u v x^{\beta} \mathrm{d}x+\int_{0}^{1} x^{\alpha+\beta}\left(u_{x} v_{x}-\frac{\mu}{x^{2}} u v\right) \mathrm{d}x=\int_{0}^{1} f v x^{\beta} \mathrm{d}x\]
for all $v\in H_{\alpha, \beta, 0}^{1}$. Therefore
$$
x^{\beta} u -\left(x^{\alpha+\beta} u_{x}\right)_x-\frac{\mu}{x^{2-\alpha-\beta}} u=x^{\beta}f \quad \text { in } \D^{\prime} \text {. }
$$
hence
$$
u-\left(x^{\alpha} u_{x}\right)_{x}-\beta x^{\alpha-1} u_{x}-\frac{\mu}{x^{2-\alpha}} u=f \text { in } \D^{\prime},
$$
thus $u\in D(\A)$ and $u+\A u=f$.\\

It follows that $-\A$ is $m$-dissipative in $L^2_\beta(0,1)$ and Corollary 2.4.10 in \cite[ p. 24]{cazenave} implies that $-\A$ is self-adjoint.\\

In \cite[Section 4]{bic} was proved that the family
\[\Psi_k(x):=\frac{2^{1/2}}{|J'_\nu(j_{\nu,k})|}x^{1/2}J_\nu(j_{\nu,k}x),\quad k\geq 1,\]
is an orthonormal basis for $L^2(0,1)$ such that 
\begin{equation}\label{zero}
\A_0\Psi_k=\tilde\lambda_k \Psi_k,\quad\text{where } \tilde\lambda_k=(j_{\nu,k})^2,\quad k\geq 1.
\end{equation}
Let $\U$ the operator given in (\ref{unitary}). Notice that $\U\Psi_k=\Fi_k$, $k\geq 1$, and the result follows by (\ref{similar}).
\end{proof}

Then $(\A,D(\A))$ generates a diagonalizable semigroup in $L^2_{\beta}(0,1)$. In fact, the last result allow us to introduce interpolation spaces for the initial data. For any $s\geq 0$, we define
\[\mathcal{H}^{s}=\mathcal{H}^{s}(0,1):=D(\mathcal{A}^{s/2})=\left\{u=\sum_{k=1}^\infty a_{k} \Phi_{k}:\|u\|_{\mathcal{H}^{s}}^{2}=\sum_{k=1}^\infty |a_{k}|^{2} \lambda_{k}^{s}<\infty\right\},\]
and we also consider the corresponding dual spaces
\[\mathcal{H}^{-s}:=\left[\mathcal{H}^{s}(0,1)\right]^{\prime}.\]
It is well known that $\mathcal{H}^{-s}$ is the dual space of $\mathcal{H}^{s}$ with respect to the pivot space $L^2_\beta(0,1)$, i.e
\[\mathcal{H}^s\hookrightarrow \mathcal{H}^0=L^2_{\beta}(0,1)=\left(L^2_{\beta}(0,1)\right)'\hookrightarrow \mathcal{H}^{-s},\quad s>0. \]
Equivalently, $\mathcal{H}^{-s}$ is the completion of $L^2_\beta(0,1)$ with respect to the norm
\[\|u\|^2_{-s}:=\sum_{k=1}^{\infty}\lambda_k^{-s}|\langle u,\Phi_k\rangle|^2.\]
It is well known that the linear mapping given by
 \[S(t)u_0=\sum_{k=1}^\infty \textrm{e}^{-\lambda_k t}a_k\Fi_k\quad\text{if}\quad u_0=\sum_{k=1}^\infty a_{k} \Phi_{k}\in \mathcal{H}^s\]
defines a semigroup $S(t)$, $t\geq 0$, in $\mathcal{H}^s $ for all $s\in\R$.\\

For $\delta\in\R$ and a smooth function $z:(0,1)\rightarrow \R$ we introduce the notion of $\delta$-generalized derivative of $z$ at $x=0$ as follows
\[\OO_\delta (z):=\lim_{x\rightarrow 0^+} x^\delta z'(x).\]

In order to get a convenient definition of weak solution for system (\ref{problem}), we multiply the equation in (\ref{problem}) by $\varphi(\tau)=S(\tau-t)z^{\tau}$, integrate by parts (formally), use the boundary-initial conditions and take the expression obtained.

\begin{definition}
Let $T>0$ and $\al,\beta,\mu\in \R$ with $0\leq\al<2$, $\al+\beta <1$, $\mu<\mu(\alpha+\beta)$. Let $f \in L^2(0,T)$ and $u_0\in \HH^{-s}$ for some $s > 0$. A weak solution of (\ref{problem}) is a function $u \in C^0([0,T];\HH^{-s})$ such that for every $\tau \in (0,T]$ and for
every $z^\tau \in \HH^s$ we have
\begin{equation}\label{weaksol}
\left\langle u(\tau), z^{\tau}\right\rangle_{\mathcal{H}^{-s}, \mathcal{H}^{s}}=\int_{0}^{\tau} f(t) \mathcal{O}_{\al+\beta+\gamma}\left(S(\tau-t) z^{\tau}\right) \mathrm{d} t+\left\langle u_{0}, S(\tau)z^\tau\right\rangle_{\mathcal{H}^{-s}, \mathcal{H}^{s}},
\end{equation}
where $\gamma=\gamma(\alpha,\beta,\mu)$ is given in (\ref{gamadef})
\end{definition}
The next result shows the existence of weak solutions for system (\ref{problem}) under suitable conditions on the parameters $\alpha,\beta,\mu,\gamma$ and $s$.
\begin{proposition}\label{continuity}
Let $T>0$ and $\al,\beta\in \R$ with $0\leq\al<2$, $\al+\beta <1$. Let $f \in L^2(0,T)$ and $u_0\in \HH^{-s}$ such that $s>\nu$, where $\nu$ is given in (\ref{Nu}). Then, formula (\ref{weaksol}) defines for each $\tau \in [0, T ]$ a unique element $u(\tau) \in \HH^{-s}$ that can be written as
\[
u(\tau)=S(\tau) u_{0}+B(\tau) f, \quad \tau \in(0, T],\]
where $B(\tau)$ is the strongly continuous family of bounded operators $B(\tau): L^{2}(0,T) \rightarrow \mathcal{H}^{-s}$ given by
 \[\left\langle B(\tau) f, z^{\tau}\right\rangle_{\mathcal{H}^{-s}, \mathcal{H}^{s}}=\int_{0}^{\tau} f(t) \mathcal{O}_{\alpha+\beta+\gamma}\left(S(\tau-t) z^{\tau}\right) \mathrm{d} t, \quad \text{for all  }z^{\tau} \in \mathcal{H}^{s} .\]
Furthermore, the unique weak solution $u$ on $[0, T]$ to (\ref{problem}) (in the sense of (\ref{weaksol})) belongs to ${C}^{0}\left([0, T] ; \mathcal{H}^{-s}\right)$ and fulfills
\[
\|u\|_{L^{\infty}\left([0, T] ; \mathcal{H}^{-s}\right)} \leq C\left(\left\|u_{0}\right\|_{\mathcal{H}^{-s}}+\|f\|_{L^{2}(0, T)}\right).
\]
\end{proposition}
\begin{proof}
Fix $\tau>0$. Let $u(\tau)\in H^{-s}$ be determined by the condition (\ref{weaksol}), hence
$$
u(\tau)-S(\tau) u_{0}=\zeta(\tau)f,
$$
where
$$
\left\langle\zeta(\tau)f, z^{\tau}\right\rangle_{\mathcal{H}^{-s}, \mathcal{H}^{s}}=\int_{0}^{\tau} f(t) \mathcal{O}_{\alpha+\beta+\gamma}\left(S(\tau-t) z^{\tau}\right) \mathrm{d} t, \quad \text{for all  } z^{\tau} \in \mathcal{H}^{s}.
$$
We claim that $\zeta(\tau)$ is a bounded operator from $L^{2}(0, T)$ into $\mathcal{H}^{-s}$: Let $z^{\tau} \in \mathcal{H}^{s}$ given by 
\begin{equation}\label{finaldata}
z^{\tau}=\sum_{k=1}^\infty a_{k} \Phi_{k},
\end{equation}
therefore
$$
S(\tau-t) z^{\tau}=\sum_{k=1}^\infty \mathrm{e}^{\lambda_{k}(t-\tau)} a_{k} \Phi_{k}, \quad \text{for all } t \in[0, \tau].
$$
By using Lemma \ref{reduce} and (\ref{limDer1}) we obtain that there exists a constant $C=C(\al,\beta,\mu)>0$ such that
\[|\mathcal{O}_{\alpha+\beta+\gamma}\left(\Phi_k\right)|\leq C |j_{\nu,k}|^{\nu+1/2},\quad k\geq 1.\]
hence (\ref{below}) implies that there exists $C=C(\al,\beta,\mu)>0$ such that
\begin{eqnarray*}
\pts{\int_{0}^{\tau}\left|\mathcal{O}_{\alpha+\beta+\gamma}\left(S(\tau-t) z^{\tau}\right)\right|^{2} \mathrm{~d} t}^{1/2}&\leq & \sum_{k=1}^\infty 
|a_k| |\OO_{\al+\beta+\gamma}(\Fi_k)| \pts{\int_0^\tau \mathrm{e}^{2\lambda_{k}(t-\tau)} \mathrm{~d} t}^{1/2} \\
&\leq& C\left\|z^{\tau}\right\|_{\mathcal{H}^{s}}\pts{\sum_{k=1}^\infty |\lambda_k|^{\nu+1/2-s}\int_0^\tau \mathrm{e}^{2\lambda_{k}(t-\tau)} \mathrm{~d} t}^{1/2}\\
&=& C\left\|z^{\tau}\right\|_{\mathcal{H}^{s}}\pts{\sum_{k=1}^\infty |\lambda_k|^{\nu-1/2-s}\pts{1-\mathrm{e}^{-2\lambda_k \tau}}}^{1/2}\\
&\leq& C\left\|z^{\tau}\right\|_{\mathcal{H}^{s}}\pts{\sum_{k=1}^\infty\frac{1}{k^{2(s-\nu+1/2)}}}^{1/2}=C\left\|z^{\tau}\right\|_{\mathcal{H}^{s}}.
\end{eqnarray*}
Therefore $\|\zeta(\tau) f\|_{\mathcal{H}^{-s}}\leq C\|f\|_{L^2(0,T)}$ for all $f\in L^2(0,T)$, $\tau\in (0,T]$.\\

Finally, we fix $f\in L^2(0,T)$ and show that the mapping $\tau\mapsto \zeta(\tau) f$ is right-continuous on $[0,T)$. Let $h>0$ small enough and $z\in \mathcal{H}^s$ given as in (\ref{finaldata}). Thus, proceeding as in the last inequalities, we have
\begin{eqnarray*}
|\left\langle\zeta(\tau+h)f-\zeta(\tau)f, z\right\rangle_{\mathcal{H}^{-s}, \mathcal{H}^{s}}|&\leq&\int_{0}^{\tau} |f(t)| |\mathcal{O}_{\alpha+\beta+\gamma}\left(\pts{S(\tau+h-t) -S(\tau-t)} z\right)| \mathrm{d} t \\
			&&+ \int_{\tau}^{\tau+h} |f(t)| |\mathcal{O}_{\alpha+\beta+\gamma}\left(S(\tau+h-t) z\right)| \mathrm{d} t\\
			&\leq& C\left\|z\right\|_{\mathcal{H}^{s}}\|f\|_{L^2(0,T)}\cts{\pts{\sum_{k=1}^\infty\frac{I(\tau,k,h)}{k^{2(s-\nu+1/2)}}}^{1/2}+\pts{\sum_{k=1}^\infty\frac{1-\mathrm{e}^{-2\lambda_k h}}{k^{2(s-\nu+1/2)}}}^{1/2}},
\end{eqnarray*}
where 
$$I(\tau,k,h)=\lambda_k\int_0^\tau\pts{\mathrm{e}^{\lambda_k(t-\tau-h)}-\mathrm{e}^{\lambda_k(t-\tau)}}^2\mathrm{~d} t\rightarrow 0\quad\text{as}\quad h\rightarrow 0^+,$$
uniformly in $k$ because $\lambda_k\geq \lambda_1$, and the result follows.
\end{proof}

\begin{remark}
In the following sections we will consider initial conditions in $L^2_\beta(0,1)$. Notice that $L^2_\beta(0,1)\subset H^{-\nu-\delta}$ for all $\delta>0$, and we can apply Proposition \ref{continuity} with $s=\nu+\delta$, $\delta>0$, then the corresponding solutions will be in $C^0([0, T ], H^{-\nu-\delta})$.
\end{remark}

\section{Upper estimate of the cost of the null controllability}
In order to solve the controllability problem we use the approach introduced by Fattorini \& Russell in \cite{Fatorini}, namely, the so-called moment method that involves the construction of a biorthogonal family $\d \lvs{\psi_k}_{k\geq 1}\subset L^2(0,T)$ to the family of exponential functions $\lvs{\mathrm{e}^{-\lambda_{k}(T-t)}}_{k\geq 1}$ on $[0, T]$, i.e that satisfies
$$
\int_{0}^{T}\psi_k(t) \mathrm{e}^{-\lambda_{l}(T-t)} dt = \delta_{kl},\quad\text{for all}\quad k,l\geq 1.
$$
As a consequence, we will get an upper bound for the cost of the null controllability.\\

Assume that each $F_k$, $k\geq 1$, is an entire function of exponential type $T/2$ such that $F_k(x)\in L^2(\R)$, and
\begin{equation}\label{krone}
F_{k}(i\lambda_{l})=\delta_{kl},\quad \text{for all}\quad k,l\geq 1.
\end{equation}
The $L^2$-version of the Paley-Wiener theorem implies that there exists $\eta_k\in L^2(\R)$ with support in $[-T/2,T/2]$ such that $F_k(z)$ is the analytic extension of the Fourier transform of $\eta_k$. Then we have that 
\begin{equation}\label{psieta}
\psi_k(t):=\mathrm{e}^{\lambda_k T/2}\eta_k(t-T/2),\quad t\in[0,T],\,k\geq1,
\end{equation}
is the family we are looking for:
\[\delta_{kl}=\mathrm{e}^{(\lambda_k-\lambda_l)T/2}F_{k}(i\lambda_{l})=\mathrm{e}^{(\lambda_k-\lambda_l)T/2}\int_{-\frac{T}{2}}^{\frac{T}{2}} \eta_{k}(t) \mathrm{e}^{\lambda_{l}t} \mathrm{d}t=\int_{0}^{T}\psi_k(t) \mathrm{e}^{-\lambda_{l}(T-t)} \mathrm{d}t\quad\text{for all}\quad k,l\geq 1.\]

Now, we proceed to construct the family $F_k$, $k\geq 1$, satisfying the aforementioned properties. First, consider the Weierstrass infinite product
\begin{equation}
	\Lambda(z):=\prod_{k=1}^{\infty}\pts{1+\dfrac{iz}{(\kappa_\al j_{\nu, k})^2}},
\end{equation}
From (\ref{asint}) we have that $j_{\nu, k}=O(k)$ for $k$ large, thus the infinite product is well-defined and converges absolutely in $\C$. Hence $\Lambda(z)$ is an entire function with simple zeros at $i(\kappa_\al j_{\nu, k})^2=i\lambda_k$, $k\geq 1$.\\

 Using \cite[Chap. XV, p. 438, eq. (3)]{Watson}, we can write
\begin{equation}\label{fLambda}
	\Lambda(z)=\Gamma(\nu+1)\pts{\dfrac{2\kappa_\al}{\sqrt{-iz}}}^\nu J_{\nu}\pts{\dfrac{\sqrt{-iz}}{\kappa_\al}}.
\end{equation}
For any $\nu>-1/2$ we have the Poisson representation for the Bessel function $J_{\nu}(z)$ (see \cite[Chap. III, p. 48, eq. (4)]{Watson}) given as follows
$$
J_{\nu}(z)=\frac{1}{\Gamma\left(\nu+\frac{1}{2}\right) \sqrt{\pi}}\left(\frac{z}{2}\right)^{\nu} \int_{-1}^{1} \mathrm{e}^{i z s}\left(1-s^{2}\right)^{\nu-\frac{1}{2}} \mathrm{d}s,\quad z\in\C,
$$
therefore
\begin{eqnarray*}
|J_{\nu}(z)|&\leq& \frac{|z|^\nu \mathrm{e}^{|\Im(z)|}}{2^\nu\Gamma\left(\nu+\frac{1}{2}\right) \sqrt{\pi}} \int_{-1}^{1} \mathrm{e}^{-(sgn(\Im(z))+s)\Im (z)}\left(1-s^{2}\right)^{\nu-\frac{1}{2}} \mathrm{d}s\\
&\leq& \frac{|z|^\nu \mathrm{e}^{|\Im(z)|}}{2^\nu\Gamma\left(\nu+\frac{1}{2}\right) \sqrt{\pi}} \int_{-1}^{1} \left(1-s^{2}\right)^{\nu-\frac{1}{2}} \mathrm{d}s= \frac{|z|^\nu \mathrm{e}^{|\Im(z)|}}{2^\nu\Gamma\left(\nu+1\right)},\quad z\in\C.
\end{eqnarray*}
The last inequality and (\ref{fLambda}) imply that
\[|\Lambda(z)|\leq \exp\pts{\frac{|\Im(\sqrt{-iz})|}{\kappa_\alpha}},\quad z\in\C.\]
In particular,
\begin{equation}\label{besst}
|\Lambda(z)|\leq\exp\pts{\frac{\sqrt{|z|}}{\kappa_\al}},\quad z\in\C,\quad |\Lambda(x)|\leq\exp\pts{\frac{\sqrt{|x|}}{\sqrt{2}\kappa_\al}},\quad x\in\R.
\end{equation}
Now, consider the entire function
\begin{equation}\label{PsiFunction}
\Psi_{k}(z):=\dfrac{\Lambda(z)}{\Lambda'(i\lambda_{k})(z-i\lambda_{k})},\quad k\geq 1.
\end{equation}
Since $\Lambda(z)$ has simple zeros at $i\lambda_k$, $k \geq 1$, we have that the family $\{\Psi_k\}$ satisfies (\ref{krone}). Unfortunately, the second estimate in (\ref{besst}) is optimal at infinity, hence $\Psi_{k}(x)$ is not in $L^2(\R)$. In order to fix this behavior, the next step is to construct a suitable ``complex multiplier". Here, we follow the approach introduced in \cite{Tucsnak}.\\

For $\theta>0$ and $a>0$, we define
$$
\sigma_{\theta}(t):=\exp\pts{-\frac{\theta}{1-t^2}},\quad t\in(-1,1),
$$
and extended by $0$ outside of $(-1, 1)$. Clearly $\sigma_{\theta}$ is analytic on $(-1,1)$. Set $C_{\theta}^{-1}:=\int_{-1}^{1}\sigma_{\theta}(t)dt$ and define
\begin{equation}\label{Hfunction}
H_{a,\theta}(z)=C_{\theta}\int_{-1}^{1}\sigma_{\theta}(t)\exp\pts{-iatz}\mathrm{d}t.
\end{equation}
Clearly $H_{a,\theta}(z)$ is a continuous function in $\C$, furthermore Morera and Fubini theorems imply that $H_{a,\theta}(z)$ is an entire function. The following result gives aditional information about $H_{a,\theta}(z)$.\\
\begin{lemma}
	The function $H_{a,\theta}$ fulfills the following inequalities
	\begin{eqnarray}
		H_{a,\theta}(ix)&\geq &\frac{\exp\pts{a|x|/\pts{2\sqrt{\theta+1}}}}{\pts{11\sqrt{\theta+1}}},\quad x\in\R,\label{Hcot1}\\
		|H_{a,\theta}(z)|     &\leq & \exp\pts{a|\Im(z)|},\quad z\in\C,\label{Hcot2}\\
		|H_{a,\theta}(x)|     &\leq & \chi_{|x|\leq 1}(x)+c\sqrt{\theta+1}\sqrt{a\theta\abs{x}}\exp\pts{3\theta/4-\sqrt{a\theta\abs{x}}}\chi_{|x|> 1}(x),\quad x\in\R,\label{Hcot3}
	\end{eqnarray}
where $c>0$ does not depend on $a$ and $\theta$.
\end{lemma}
We refer to \cite[pp. 85--86]{Tucsnak} for the details.\\

Finally, for $k\geq 1$ consider the entire function $F_{k}$ given as
\begin{equation}\label{Ffunction}
	F_{k}(z):=\Psi_{k}(z)\dfrac{H_{a,\theta}(z)}{H_{a,\theta}(i\lambda_{k})},\quad z\in\C.
\end{equation}
 
For $\delta\in(0,1)$ we set
\begin{equation}\label{aConst}
	a:=\frac{T(1-\delta)}{2}>0,\quad \text{and}\quad \theta:=\dfrac{(1+\delta)^2}{\kappa_\al^2 T\pts{1-\delta}}>0.
\end{equation}

The computations in the next result justify the choice of these parameters.
\begin{lemma}
	For each $k\geq 1$ the function $F_{k}(z)$ satisfies the following properties:\\
	i) $F_{k}$ is of exponential type $T/2$,\\
	ii) $F_{k}\in L^1(\R)\cap L^2(\R)$,\\
	iii) $F_k$ satisfies (\ref{krone}).\\
	iv) Furthermore, there exists a constant $c>0$, independent of $T,\alpha$ and $\delta$, such that
	\begin{equation}\label{Fbound}
		\left\|F_{k}\right\|_{L^{1}(\R)} \leq \frac{C(T, \alpha,\delta)}{\lambda_{k}\left|\Lambda^{\prime}\left(i \lambda_{k}\right)\right|}  \exp\pts{-\frac{a\lambda_k}{2\sqrt{\theta+1}}},		
	\end{equation} 
	where
	\begin{equation}\label{upper}
		C(T, \alpha,\delta)=c\sqrt{\theta+1}\cts{\exp\pts{{\frac{1}{\sqrt{2}\kappa_\alpha}}}+\sqrt{\theta+1}\frac{\kappa_\alpha^2}{\delta^3}\exp\pts{\frac{3 \theta}{4}}}. 
		\end{equation}
\end{lemma}
\begin{proof}
	By using (\ref{besst}), (\ref{Hcot2}), (\ref{Ffunction}) and (\ref{aConst}) we get that $F_{k}$ is of exponential type $T/2$. Moreover, by (\ref{PsiFunction}), (\ref{Hcot1}) and (\ref{Ffunction}), $F_{k}$ fulfills (\ref{krone}).\\
	
	Now we use (\ref{besst}), (\ref{Hcot1}), (\ref{Hcot3}), (\ref{Ffunction}), and (\ref{aConst}) to get
	\begin{eqnarray*}
	\left|F_{k}(x)\right| &\leq& c\exp\pts{-\frac{a\lambda_k}{2\sqrt{\theta+1}}}\frac{\sqrt{\theta+1}}{|\Lambda^{\prime}\left(i \lambda_{k}\right)||x^2+ \lambda_{k}^2|^{1/2}} |H_{a,\theta}(x)|\exp\pts{ \frac{|x|^{1/2}}{\sqrt{2}\kappa_\alpha}}
						\\
				&\leq& c \exp\pts{-\frac{a\lambda_k}{2\sqrt{\theta+1}}}\frac{\sqrt{\theta+1}}{\lambda_{k}|\Lambda^{\prime}\left(i \lambda_{k}\right)|}\cts{e^{\frac{1}{\sqrt{2}\kappa_\alpha}}\chi_{|x|\leq 1}(x)+\sqrt{\theta+1}\sqrt{a\theta\abs{x}}\exp\pts{\frac{3\theta}{4}- \frac{\delta|x|^{1/2}}{\sqrt{2}\kappa_\alpha}}\chi_{|x|> 1}(x)}.
	\end{eqnarray*}	
Since the function in the right hand side is rapidly decreasing in $\R$, we have $F_{k}\in L^1(\R)\cap L^2(\R)$. Finally, the change of variable $y=(\kappa_\alpha)^{-1}\delta|x|^{1/2}/\sqrt{2}$ implies (\ref{Fbound}).
\end{proof}

Since $\eta_k, F_k\in L^1(\R)$, the inverse Fourier theorem yields 
\[\eta_k(t)=\frac{1}{2\pi}\int_{\R}\mathrm{e}^{it\tau}F_k(\tau)\mathrm{d}\tau,\quad t\in\R, k\geq 1,\]
hence (\ref{psieta}) implies that $\psi_k\in C([0,T])$ and from (\ref{Fbound}) we have
\begin{equation}\label{psiL1}
\|\psi_k\|_{\infty}\leq \frac{C(T, \alpha,\delta)}{\lambda_{k}\left|\Lambda^{\prime}\left(i \lambda_{k}\right)\right|}  \exp\pts{\frac{T\lambda_k}{2}-\frac{a\lambda_k}{2\sqrt{\theta+1}}},\quad k\geq 1.
\end{equation}

Now, we are ready to prove the null controllability of system (\ref{problem}). Let $u_{0}\in L^{2}_\beta(0,1)$. Then consider its Fourier series with respect to the orthonormal basis $\{\Phi_{k}\}_{k\geq 1}$,
 \begin{equation}\label{uoSerie}
 	u_{0}(x)=\sum_{k=1}^{\infty} a_{k} \Phi_{k}(x).
 \end{equation}
 We set
 \begin{equation}\label{fserie}
 	f(t):=-\sum_{k=1}^{\infty}\frac{a_{k} \mathrm{e}^{-\lambda_{k} T}}{\mathcal{O}_{\alpha+\beta+\gamma}\left(\Phi_{k}\right)} \psi_{k}(t)
 \end{equation}
provided that $\alpha+\beta<1$ and $\gamma$ given in (\ref{gamadef}). Since $\{\psi_k\}$ is biorthogonal to $\{\mathrm{e}^{-\lambda_k(T-t)}\}$ we have
$$
\int_{0}^{T} f(t) \mathcal{O}_{\alpha+\beta+\gamma}\left(\Phi_{k}\right) \mathrm{e}^{-\lambda_{k}(T-t)} \mathrm{d} t=-a_{k} \mathrm{e}^{-\lambda_{k} T}=-\left\langle u_{0}, \mathrm{e}^{-\lambda_{k} T}\Phi_{k}\right\rangle_\beta=-\left\langle u_{0}, \mathrm{e}^{-\lambda_{k} T}\Phi_{k}\right\rangle_{\mathcal{H}^{-s}, \mathcal{H}^{s}}.
$$
Let $u\in C([0,T];H^{-s})$ that satisfies (\ref{weaksol}) for all $\tau\in (0,T]$, $z^\tau\in H^s$. In particular, for $\tau=T$ we take $z^T=\Phi_k$, $k\geq 1$, then the last equality implies that
$$
\left\langle u(\cdot, T), \Phi_{k}\right\rangle_{\mathcal{H}^{-s}, \mathcal{H}^{s}}=0\quad \text{for all}\quad k \geq 1,
$$
hence $u(\cdot, T)=0$.\\

It just remains to estimate the norm of the control $f$. From (\ref{psiL1}) and (\ref{fserie})  we get
\begin{equation}\label{finty}
	\|f\|_{\infty} \leq  C(T, \alpha,\delta)\sum_{k=1}^{\infty} \frac{\left|a_{k}\right|}{\left|\mathcal{O}_{\alpha+\beta+\gamma}\left(\Phi_{k}\right)\right|} \frac{1}{\lambda_{k}\left|\Lambda^{\prime}\left(i \lambda_{k}\right)\right|} \exp\pts{-\frac{T\lambda_k}{2}-\frac{a\lambda_k}{2\sqrt{\theta+1}}}.
\end{equation}

From (\ref{fLambda}), we obtain
\begin{equation*}\label{estimate1}
	\left|\Lambda^{\prime}\left(i \lambda_{k}\right)\right|=\Gamma(\nu+1)\frac{2^{\nu}}{|j_{\nu, k}|^{\nu}} \frac{1}{2\kappa_\al^{2} j_{\nu, k}} |J_{\nu}^{\prime}\left(j_{\nu, k}\right)|, \quad k\geq 1.
\end{equation*}
and by using (\ref{lambdak}) and (\ref{limDer1}) we get
\begin{equation*}\label{estimate2}
	\d\left|\mathcal{O}_{\alpha+\beta+\gamma}\left(\Phi_{k}\right)\lambda_{k}\Lambda^{\prime}\left(i \lambda_{k}\right)\right|= 2^{-1/2}\sqrt{\kappa_\al}(\sqrt{\mu(\alpha+\beta)}+\kappa_\al\nu) j_{\nu, k}. 
\end{equation*}
From (\ref{finty}), the last equality and using that $\lambda_k\geq \lambda_1$, it follows that 
\begin{equation*}
	\|f\|_{\infty} \leq \frac{C(T, \alpha,\delta)\kappa_\al^{-1/2}}{\sqrt{\mu(\alpha+\beta)}+\kappa_\al\nu} \exp\pts{-\frac{T\lambda_1}{2}-\frac{a\lambda_1}{2\sqrt{\theta+1}}}\sum_{k=1}^{\infty} \frac{|a_{k}| }{j_{\nu, k}}.
	\end{equation*}
By using the Cauchy-Schwarz inequality, the fact that $j_{\nu,k}\geq (k-1/4)\pi$ (by (\ref{below}))  and (\ref{uoSerie}), we obtain that
\begin{eqnarray*}
	\|f\|_{\infty} &\leq& \frac{C(T, \alpha,\delta)\kappa_\al^{-1/2}}{\sqrt{\mu(\alpha+\beta)}+\kappa_\al\nu} \exp\pts{-\frac{T\lambda_1}{2}-\frac{a\lambda_1}{2\sqrt{\theta+1}}}\left(\sum_{k=1}^{\infty} |a_{k}|^{2}\right)^{\frac{1}{2}}\\
	&=&\frac{C(T, \alpha,\delta)\kappa_\al^{-1/2}}{\sqrt{\mu(\alpha+\beta)}+\sqrt{\mu(\alpha+\beta)-\mu)}} \exp\pts{-\frac{T\lambda_1}{2}-\frac{a\lambda_1}{2\sqrt{\theta+1}}}\left\|u_{0}\right\|_{\beta}.
\end{eqnarray*}
Using the expression of $a,\theta$ given in (\ref{aConst}) and the facts $\theta>0$, $\delta\in(0,1)$, and $0<\kappa_\alpha\leq 1$, we get that
\[\theta\leq \frac{4}{(1-\delta)\kappa_\alpha^2 T },\quad \sqrt{\theta+1}\leq \frac{2(1+T)^{1/2}}{(1-\delta)^{1/2}\kappa_\alpha T^{1/2}},\quad \sqrt{\theta+1}\leq \theta+1,\]
therefore
\[\frac{a}{\sqrt{\theta+1}}\geq \frac{\kappa_\alpha(1-\delta)^{3/2}T^{3/2}}{4(1+T)^{1/2}}, \quad C(T, \alpha,\delta)\leq c\pts{1+\frac{1}{(1-\delta)\kappa_\alpha^2 T}}\cts{\exp\pts{\frac{1}{\sqrt{2}\kappa_\alpha}}+\frac{1}{\delta^3}\exp\pts{\frac{3}{(1-\delta)\kappa_\alpha^2 T}}}, \]
and by using the definition of $\lambda_{1}$ the result follows.

\section{Lower estimate of the cost of the null controllability}
In this section we get a lower estimate of the cost $\mathcal{K}=\mathcal{K}(T,\alpha,\beta,\mu)$.
We set
\begin{equation}\label{first} 
u_0(x):= \frac{\abs{J'_\nu\pts{j_{\nu,1}}}}{(2\kappa_\al)^{1/2}}\Fi_1(x),\,x\in(0,1),
\quad \text{hence}\quad\|u_0\|^2_\beta=\frac{\abs{J'_\nu\pts{j_{\nu,1}}}^2}{2\kappa_\al}.
\end{equation}
For $\varepsilon>0$ small enough, there exists $f\in U(\al,\beta,\mu,T,u_0)$ such that
\begin{equation}\label{inicost}
u(\cdot,T)\equiv 0,\quad \text{and}\quad \|f\|_{L^2(0,T)}\leq (\mathcal{K}+\varepsilon)\|u_0\|_\beta.
\end{equation}
Then, in (\ref{weaksol}) we set $\tau=T$ and  take $z^\tau=\Fi_k$, $k\geq 1$,  to obtain
\begin{eqnarray*}
\mathrm{e}^{-\lambda_k T}\left\langle u_{0},\Fi_k\right\rangle_\beta=\left\langle u_{0}, S(T)\Fi_k\right\rangle_{\mathcal{H}^{-s}, \mathcal{H}^{s}}&=&
-\int_{0}^{T} f(t) \mathcal{O}_{\al+\beta+\gamma}\left(S(T-t) \Fi_k\right) \mathrm{d} t\\
&=&
-\mathrm{e}^{-\lambda_k T}\mathcal{O}_{\al+\beta+\gamma}\left(\Fi_k\right)\int_{0}^{T} f(t) \mathrm{e}^{\lambda_k t} \mathrm{d} t,
\end{eqnarray*}
from (\ref{first}) and (\ref{limDer1}) it follows that
\begin{equation}\label{orto}
\int_{0}^{T} f(t) \mathrm{e}^{\lambda_k t} \mathrm{d} t= -\frac{2^{\nu}\Gamma(\nu+1)\abs{J'_{\nu}\pts{j_{\nu, 1}}}^2}{2\kappa_\al\pts{\sqrt{\mu(\alpha+\beta)}+\kappa_\al \nu}\pts{j_{\nu,1}}^{\nu}}\delta_{1,k},\quad k\geq 1.
\end{equation}

Now consider the function $v: \mathbb{C} \rightarrow \mathbb{C}$ given by
\begin{equation*}
v(s):=\int_{-T / 2}^{T / 2} f\left(t+\frac{T}{2}\right) \mathrm{e}^{-i s t} \mathrm{~d} t, \quad s \in \mathbb{C} .
\end{equation*}
Fubini and Morera theorems imply that $v(s)$ is an entire function. Moreover, (\ref{orto}) implies that
\[v(i\lambda_k)=0\quad\text{for all }k\geq 2,\quad \text{and}\quad v(i\lambda_1)=-\frac{2^{\nu}\Gamma(\nu+1)\abs{J'_{\nu}\pts{j_{\nu, 1}}}^2}{2\kappa_\al\pts{\sqrt{\mu(\alpha+\beta)}+\kappa_\al \nu}\pts{j_{\nu,1}}^{\nu}}\mathrm{e}^{-\lambda_1 T/2}.\]
We also have that
\begin{equation}\label{uve}
|v(s)| \leq \mathrm{e}^{T|\Im(s)|/2}\int_{0}^{T}|f(t)| \mathrm{d} t \leq (\mathcal{K}+\varepsilon) T^{1/2}\mathrm{e}^{T|\Im(s)|/2} \left\|u_{0}\right\|_{\beta}.
\end{equation}
Consider the entire function $F(z)$ given by
\begin{equation}\label{entire}
F(s):=v\left(s-i \delta\right), \quad s \in \mathbb{C},
\end{equation}
for some $\delta>0$ that will be chosen later on. Clearly, 
\[ F\left(a_{k}\right)=0, \quad k\geq 2, \quad \text {where} \quad a_{k}:=i\left(\lambda_k+\delta\right),\quad k\geq 1,\quad\text{and}\]
\begin{equation}\label{ena1}
F\left(a_{1}\right)=-\frac{2^{\nu}\Gamma(\nu+1)\abs{J'_{\nu}\pts{j_{\nu, 1}}}^2}{2\kappa_\al\pts{\sqrt{\mu(\alpha+\beta)}+\kappa_\al \nu}\pts{j_{\nu,1}}^{\nu}}\mathrm{e}^{-\lambda_1 T/2}.
\end{equation}
From (\ref{first}), (\ref{uve}) and (\ref{entire}) we obtain
\begin{equation}\label{logF}
\log |F(s)|\leq \frac{T}{2}|\Im(s)-\delta|+\log\pts{(\mathcal{K}+\varepsilon) T^{1 / 2}\frac{\abs{J'_\nu\pts{j_{\nu,1}}}}{\pts{2\kappa_\al}^{1/2}}},\quad s\in\mathbb{C}.
\end{equation}
We recall the following representation theorem, see \cite[p. 56]{koosis}.
\begin{theorem} Let $g(z)$ be an entire function of exponential type and assume that
$$
\int_{-\infty}^{\infty} \frac{\log ^{+}|g(x)|}{1+x^{2}} \mathrm{d}x<\infty.
$$
Let $\left\{b_{\ell}\right\}_{\ell \geq 1}$ be the set of zeros of $g(z)$ in the upper half plane $\Im(z)>0$ (each zero being repeated as many times as its multiplicity). Then,
$$
\log |g(z)|=A \Im(z)+\sum_{\ell=1}^{\infty} \log \left|\frac{z-b_{\ell}}{z-\bar{b}_{\ell}}\right|+\frac{\Im(z)}{\pi} \int_{-\infty}^{\infty} \frac{\log |g(s)|}{|s-z|^{2}} \mathrm{d}s,\quad\Im(z)>0,
$$
where
$$
A=\limsup _{y \rightarrow\infty} \frac{\log |g(i y)|}{y} .
$$
\end{theorem}
We apply the last result to the function $F(z)$ given in (\ref{entire}). In this case, (\ref{uve}) implies that $A\leq T/2$. Also notice that $\Im\left(a_{k}\right)>0$, $k\geq 1$, to get
\begin{equation}\label{aprep}
\log \left|F\left(a_{1}\right)\right|\leq\left(\lambda_1+\delta\right)\frac{T}{2}+\sum_{k=2}^{\infty} \log \left|\frac{a_{1}-a_{k}}{a_{1}-\bar{a}_{k}}\right|+\frac{\Im\left(a_{1}\right)}{\pi} \int_{-\infty}^{\infty} \frac{\log |F(s)|}{\left|s-a_{1}\right|^{2}} \mathrm{~d}s.
\end{equation}
By using the definition of the constants $a_k$'s we have
\begin{eqnarray}
\sum_{k=2}^{\infty} \log \left|\frac{a_{1}-a_{k}}{a_{1}-\bar{a}_{k}}\right|&=&\sum_{k=2}^{\infty} \log \left(\frac{\left( j_{\nu, k}\right)^{2}-\left( j_{\nu, 1}\right)^{2}}{2 \delta / \kappa_\alpha^{2}+\left( j_{\nu, 1}\right)^{2}+\left(j_{\nu, k}\right)^{2}}\right)\notag\\
&\leq& \sum_{k=2}^{\infty} \frac{1}{j_{\nu, k+1}-j_{\nu, k}} \int_{j_{\nu, k}}^{j_{\nu, k+1}} \log \left(\frac{ x^{2}}{2 \delta / \kappa_\alpha^{2}+ x^{2}}\right) \mathrm{d} x \label{apoyo1} \\ 
&\leq& \frac{1}{\pi} \int_{j_{\nu, 2}}^{\infty} \log \left(\frac{ x^{2}}{2 \delta / \kappa_\alpha^{2}+x^{2}}\right) \mathrm{d} x,\notag\\
&=& -\frac{j_{\nu, 2}}{\pi} \log \left(\frac{1}{1+2 \delta /\left(\kappa_\alpha j_{\nu, 2}\right)^{2}}\right)- \frac{2\sqrt{2 \delta}}{\pi\kappa_\alpha}\left(\frac{\pi}{2}-\tan ^{-1}\left(\kappa_\alpha j_{\nu, 2} / \sqrt{2 \delta}\right)\right) ,\notag
\end{eqnarray}
where we have used Lemma \ref{consec} and made the change of variables
$$
\tau=\frac{ \kappa_\alpha}{\sqrt{2 \delta}} x.
$$
From (\ref{logF}) we get the estimate
\begin{equation}\label{apoyo2}
\frac{\Im\left(a_{1}\right)}{\pi} \int_{-\infty}^{\infty} \frac{\log |F(s)|}{\left|s-a_{1}\right|^{2}} \mathrm{~d} s \leq \frac{T \delta}{2}+\log \left((\mathcal{K}+\varepsilon) T^{1 / 2} \frac{\left|J_{\nu}^{\prime}\left(j_{\nu, 1}\right)\right|}{\pts{2 \kappa_\alpha}^{1/2}}\right).
\end{equation}
From (\ref{ena1}), (\ref{aprep}), (\ref{apoyo1}) and (\ref{apoyo2}) we have
\[ \frac{2\sqrt{2 \delta}}{\pi\kappa_\alpha}\tan ^{-1}\left(\frac{\sqrt{2 \delta}}{\kappa_\alpha j_{\nu, 2}}\right) -\frac{j_{\nu, 2}}{\pi} \log \left(1+\frac{2 \delta}{ \left(\kappa_\alpha j_{\nu, 2}\right)^{2}}\right) 
-\pts{\lambda_1+\delta}T \leq \log(\mathcal{K}+\varepsilon)+\log h(\alpha,\beta,\mu, T),\]
where
\[h(\alpha,\beta,\mu, T)=\frac{\pts{{2T \kappa_\alpha}}^{1/2}\pts{\sqrt{\mu(\alpha+\beta)}+\kappa_\al \nu}\left(j_{\nu, 1}\right)^{\nu}}{2^{\nu} \Gamma(\nu+1) \left|J_{\nu}^{\prime}\left(j_{\nu, 1}\right)\right|}=\frac{\pts{{2T \kappa_\alpha}}^{1/2}\pts{\sqrt{\mu(\alpha+\beta)}+\sqrt{\mu(\alpha+\beta)-\mu}}\left(j_{\nu, 1}\right)^{\nu}}{2^{\nu} \Gamma(\nu+1) \left|J_{\nu}^{\prime}\left(j_{\nu, 1}\right)\right|}.\]
The results follows by taking 
\[\delta=\frac{\kappa_\alpha^2 \pts{j_{\nu,2}}^2}{2}, \quad \text{and then letting}\quad \varepsilon\rightarrow 0^+.\]

\begin{remark}\label{error}
a) In this item we take the definitions from \cite{bic}. In the middle of page 207 in \cite{bic}, the author uses implicitly that
\begin{equation}\label{biccerror}
\left.\left(x^{\alpha} \frac{p(x)}{p(0)}\right)^{\prime} \Phi_{k}(x)\right|_{0} ^{1}=0,
\end{equation}
where $\alpha=(1-\sqrt{1-4\mu})/2$, $p(x)=1-x^{1-2\alpha}$, $\Phi_{k}(x)=\sqrt{2}|J'_\nu(j_{\nu,k})|^{-1}x^{1/2}J_\nu(j_{\nu,k}x)$, $\nu=1/2-\alpha$. The expresion in (\ref{biccerror}) does vanish at $x=1$ but does not at $x=0$, by (\ref{asincero}). Thus the proof of Theorem 2.1 in \cite{bic} is not complete.\\

b) Here we take the definitions from \cite{bic2}. In the middle of page 525 in \cite{bic2}, the authors use implicitly that
\begin{equation}\label{biccerror2}
\left.\left(x^{\gamma} \frac{p(x)}{p(0)}\right)^{\prime} x^\alpha\Phi_{k}(x)\right|_{0} ^{1}=0,
\end{equation}
where $\gamma=\sqrt{\mu(\al)}-\sqrt{\mu(\al)-\mu}$, $p(x)=1-x^{q}$, $q=2\sqrt{\mu(\al)-\mu}$, $\Phi_{k}(x)=\sqrt{2-\al}|J'_\nu(j_{\nu,k})|^{-1}x^{\frac{1-\al}{2}}J_\nu(j_{\nu,k}x^{\frac{2-\al}{2}})$, $\nu=2\sqrt{\mu(\al)-\mu}/(2-\al)$, $0\leq\al<1$. The expression in (\ref{biccerror2}) does vanish at $x=1$ but does not at $x=0$, by (\ref{asincero}). Thus the proof of Theorem 2.2 in \cite{bic2} is not complete.
\end{remark}

\appendix
\section{Bessel functions}
We introduce the Bessel function of the first kind $J_{\nu}$ as follows
\begin{equation}\label{bessel}
J_{\nu}(x)=\sum_{m \geq 0} \frac{(-1)^{m}}{m ! \Gamma(m+\nu+1)}\left(\frac{x}{2}\right)^{2 m+\nu}, \quad x \geq 0,
\end{equation}
where $\Gamma(\cdot)$ is the Gamma function. In particular, for $\nu>-1$ and $0<x \leq \sqrt{\nu+1}$, from (\ref{bessel}) we have (see \cite[9.1.7, p. 360]{abram})
\begin{equation}\label{asincero}
J_{\nu}(x) \sim \frac{1}{\Gamma(\nu+1)}\left(\frac{x}{2}\right)^{\nu} \quad \text { as } \quad x \rightarrow 0^{+} .
\end{equation}
Bessel functions of the first kind satisfy the recurrence formula $([1], 9.1 .27)$ :
\begin{equation}\label{recur}
x J_{\nu}^{\prime}(x)-\nu J_{\nu}(x)=-x J_{\nu+1}(x).
\end{equation}
Recall the asymptotic behavior of the Bessel function $J_{\nu}$ for large $x$, see \cite[Lem. 7.2, p. 129]{komo}.
\begin{lem}\label{asimxinf}
For any $\nu \in \mathbb{R}$
$$
J_{\nu}(x)=\sqrt{\frac{2}{\pi x}}\left\{\cos \left(x-\frac{\nu \pi}{2}-\frac{\pi}{4}\right)+\mathcal{O}\left(\frac{1}{x}\right)\right\} \quad \text { as } \quad x \rightarrow \infty.
$$
\end{lem}

For $\nu >-1$ the Bessel function $J_{\nu}$ has an infinite number of real zeros $0<j_{\nu, 1}<j_{\nu, 2}<\ldots$, all of which are simple, with the possible exception of $x=0$. In \cite[Proposition 7.8]{komo} we can find the next information about the location of the zeros of the Bessel functions $J_{\nu}$:
\begin{lem}\label{consec}Let $\nu \geq 0$.\\
1. The difference sequence $\left(j_{\nu, k+1}-j_{\nu, k}\right)_{k}$ converges to $\pi$ as $k \rightarrow\infty$.\\
2. The sequence $\left(j_{\nu, k+1}-j_{\nu, k}\right)_{k}$ is strictly decreasing if $|\nu|>\frac{1}{2}$, strictly increasing if $|\nu|<\frac{1}{2}$, and constant if $|\nu|=\frac{1}{2}$.\\
\end{lem}

For $\nu \geq 0$ fixed, we consider the next asymptotic expansion of the zeros of the Bessel function $J_{\nu}$, see\cite[Section 15.53]{Watson},
\begin{equation}\label{asint}
j_{\nu, k}=\left(k+\frac{\nu}{2}-\frac{1}{4}\right) \pi-\frac{4 \nu^{2}-1}{8\left(k+\frac{\nu}{2}-\frac{1}{4}\right) \pi}+O\left(\frac{1}{k^{3}}\right), \quad \text { as } k \rightarrow\infty.
\end{equation}

In particular we have
\begin{equation}\label{below}
\begin{aligned}
&j_{\nu, k} \geq\left(k-\frac{1}{4}\right) \pi \quad \text { for } \nu \in\left[0, 1/2\right], \\
&j_{\nu, k} \geq\left(k-\frac{1}{8}\right) \pi \quad \text { for } \nu \in\left[1/2,\infty\right].
\end{aligned}
\end{equation}

\begin{lem}\label{reduce} For any $\nu \geq 0$ and any $k\geq 1$ we have
$$
\sqrt{j_{\nu, k}}\left|J_{\nu}^{\prime}\left(j_{\nu, k}\right)\right|=\sqrt{\frac{2}{\pi}}+O\left(\frac{1}{j_{\nu, k}}\right)\quad \text{as}\quad k \rightarrow \infty.
$$
\end{lem}
\begin{lem} Let $\gamma=\gamma(\alpha,\beta,\mu)$ and $\nu=\nu(\alpha,\beta,\mu)$ given in (\ref{gamadef}) and (\ref{Nu}) respectively, then the $\al+\beta+\gamma$-generalized derivative of $\Phi_k$ at $x=0$ is finite for all $k\geq 1$, and
	\begin{equation}\label{limDer1}
	\OO_{\al+\beta+\gamma}(\Fi_k)= \frac{(2\kappa_\al)^{1/2}\pts{\sqrt{\mu(\alpha+\beta)}+\kappa_\al \nu}\pts{j_{\nu,k}}^{\nu}}{2^{\nu}\Gamma(\nu+1)\abs{J'_{\nu}\pts{j_{\nu, k}}}}, \quad k\geq1.
	\end{equation}
\end{lem}
The proof of this result follows by using  (\ref{asincero}) and the recurrence formula (\ref{recur}).

\end{document}